\documentclass[leqno,a4paper,12pt]{article}

\newtheorem{Thm}{\indent Theorem}[section]
\newtheorem{Prop}[Thm]{\indent Proposition}
\newtheorem{Lem}[Thm]{\indent Lemma}

\newtheorem{Cor}[Thm]{\indent Corollary}
\newtheorem{Var}[Thm]{\indent Variant}

\newtheorem{Rem}[Thm]{\indent Remark}

\newenvironment{myequation}
               {\addtocounter{Thm}{1}\begin{equation}}
               {\end{equation}}

\newcommand{\myheading}[1]{\medskip{\bf{#1}}\quad}

\def\qed{{\hskip0pt\unskip\unskip\nobreak\hfil\penalty50
          \hskip1em\hbox{}\nobreak\hfil
          {\bf q.e.d.}%
          \parfillskip=0pt\finalhyphendemerits=0
          \par}\medskip}

\newenvironment{Proof}
               {{\it Proof.}\quad}
               {\qed}

\newcommand{\Prime}{\kern3\fontdimen1\font$'$\kern-7\fontdimen1\font}

\long\def\forget#1{}

\long\def\beginSIDEREMARK#1\endSIDEREMARK
    {{\par\bigskip\advance\leftskip by 2cm
                  \advance\rightskip by -2cm\noindent
      {\bf Our own side remark:} #1
      \par\bigskip\noindent}}

\long\def\beginFORGET#1\endFORGET{#1}
\long\def\beginFORGET#1\endFORGET{}

\def\?{\ ???\ \immediate\write16{}%
\immediate\write16{Warning: There was still a question mark . . . }%
\immediate\write16{}}

\usepackage{exscale}
\usepackage{amssymb}

\newcommand{\BA}{{\mathbb{A}}}

\newcommand{\BC}{{\mathbb{C}}}

\newcommand{\BF}{{\mathbb{F}}}
\newcommand{\BG}{{\mathbb{G}}}

\newcommand{\BQ}{{\mathbb{Q}}}
\newcommand{\BR}{{\mathbb{R}}}
\newcommand{\BS}{{\mathbb{S}}}

\newcommand{\BV}{{\mathbb{V}}}
\newcommand{\BW}{{\mathbb{W}}}

\newcommand{\BZ}{{\mathbb{Z}}}

\newcommand{\FS}{{\mathfrak{S}}}

\newcommand{\FU}{{\mathfrak{U}}}
\newcommand{\FV}{{\mathfrak{V}}}

\newcommand{\FX}{{\mathfrak{X}}}

\newcommand{\CC}{{\cal C}}

\newcommand{\CF}{{\cal F}}

\newcommand{\CH}{{\cal H}}

\newcommand{\CM}{{\cal M}}

\newcommand{\CR}{{\cal R}}
\newcommand{\CS}{{\cal S}}
\newcommand{\CT}{{\cal T}}

\newcommand{\CY}{{\cal Y}}

\newcommand{\Spec}{\mathop{{\bf Spec}}\nolimits}

\newcommand{\Gal}{\mathop{\rm Gal}\nolimits}
\newcommand{\GL}{{\rm GL}}

\newcommand{\Gr}{{\rm Gr}}
\newcommand{\iint}{\mathop{{\rm int}}\nolimits}
\newcommand{\image}{\mathop{{\rm Im}}\nolimits}
\newcommand{\imm}{\mathop{{\rm im}}\nolimits}

\newcommand{\Mod}{\mathop{\rm Mod}\nolimits}
\newcommand{\Sym}{\mathop{\rm Sym}\nolimits}
\newcommand{\Lie}{\mathop{\rm Lie}\nolimits}

\newcommand{\Hom}{\mathop{\rm Hom}\nolimits}

\newcommand{\Cent}{\mathop{\rm Cent}\nolimits}

\newcommand{\Stab}{\mathop{\rm Stab}\nolimits}

\newcommand{\Res}{\mathop{\rm Res}\nolimits}
\newcommand{\Ind}{\mathop{\rm Ind}\nolimits}

\newcommand{\Tor}{\mathop{\rm Tor}\nolimits}

\newcommand{\OFU}{\overline{\FU}}

\newcommand{\pr}{\mathop{\rm pr}\nolimits}

\def\tei{\, | \,}
\def\halb{\frac{1}{2}}

\def\phi{\varphi}
\def\epsilon{\varepsilon}

\def\id{{\rm id}}

\newbox\mybox
\def\arrover#1{\mathrel{
       \setbox\mybox=\hbox spread 1.4em{\hfil$\scriptstyle#1$\hfil}
       \vbox{\offinterlineskip\copy\mybox
             \hbox to\wd\mybox{\rightarrowfill}}}}
\def\larrover#1{\mathrel{
       \setbox\mybox=\hbox spread 1.4em{\hfil$\scriptstyle#1$\hfil}
       \vbox{\offinterlineskip\copy\mybox
             \hbox to\wd\mybox{\leftarrowfill}}}}

\def\ontoover#1{\mathrel{
       \setbox\mybox=\hbox spread 1.4em{\hfil$\scriptstyle#1$\hfil}
       \vbox{\offinterlineskip\copy\mybox
             \hbox to\wd\mybox{\rightarrowfill\hskip-2.8mm
                               $\rightarrow$}}}}
\def\leftontoover#1{\mathrel{
       \setbox\mybox=\hbox spread 1.4em{\hfil$\scriptstyle#1$\hfil}
       \vbox{\offinterlineskip\copy\mybox
             \hbox to\wd\mybox{$\leftarrow$\hskip-2.8mm
                               \leftarrowfill}}}}
\def\longto{\longrightarrow}
\def\into{\hookrightarrow}

\def\longonto{\ontoover{\ }}
\def\isoto{\arrover{\sim}}

\usepackage[curve,matrix,arrow,cmtip]{xy}
\NoComputerModernTips

\def\myxymessage{\def\messagetext
   {Here an xy-pic diagram was omitted to speed up compilation . . . }
   \immediate\write16{\messagetext}
   \hbox{\bf \messagetext}}
\def\filxymatrix#1{\myxymessage}
\def\filxyarray#1{\myxymessage}

\newdir^{ (}{{}*!/-3pt/\dir^{(}}    
\newdir_{ (}{{}*!/-3pt/\dir_{(}}    
\newdir^{ )}{{}*!/+3pt/\dir^{)}}    
\newdir_{ )}{{}*!/+3pt/\dir_{)}}    

\def\rscript#1{\hbox to 0pt{$\scriptstyle#1$\hss}}

\newcommand{\red}{{\rm red}}

\newcommand{\ad}{{\rm ad}}
\newcommand{\topp}{\rm{top}}

\newcommand{\Mps}{\mathop{M^{\pi_{[\sigma]} (K_1)}}\nolimits}
\newcommand{\MpC}{\mathop{M_{\BC}^{\pi_{[\sigma]} (K_1)}}\nolimits}
\newcommand{\Pes}{\mathop{P_{1,[\sigma]}}\nolimits}
\newcommand{\Ses}{\mathop{\FS_{1,[\sigma]}}\nolimits}
\newcommand{\Xes}{\mathop{\FX_{1,[\sigma]}}\nolimits}
\newcommand{\Mpt}{\mathop{M^{\pi_{[\tau]} (K_1)}}\nolimits}

\newcommand{\Mut}{\mathop{\mu_{K,\topp}}\nolimits}

\newcommand{\Pk}{\mathop{\Perv_F M^K _{\BC}}\nolimits}

\newcommand{\PK}{\mathop{\Perv_F \MpC}\nolimits}

\newcommand{\SP}{\mathop{Sp}\nolimits}
\newcommand{\Ss}{\mathop{\SP_\sigma}\nolimits}

\newcommand{\Rep}{\mathop{\bf Rep}\nolimits}

\newcommand{\Et}{\mathop{\bf Et}\nolimits}
\newcommand{\Loc}{\mathop{\bf Loc}\nolimits}
\newcommand{\Perv}{\mathop{\bf Perv}\nolimits}
\newcommand{\MHM}{\mathop{\bf MHM}\nolimits}
\newcommand{\VMHS}{\mathop{\bf Var}\nolimits}

\newcommand{\con}{{\rm c}}

\let\oldbullet\bullet
\def\bullet{{\mathchoice{\oldbullet}%
                        {\oldbullet}%
                        {\scriptscriptstyle\oldbullet}%
                        {\oldbullet}}}

\begin{document}  

%

\hfuzz=3pt
\overfullrule=10pt                   


\setlength{\abovedisplayskip}{6.0pt plus 3.0pt}
\setlength{\belowdisplayskip}{6.0pt plus 3.0pt}     
\setlength{\abovedisplayshortskip}{6.0pt plus 3.0pt}
\setlength{\belowdisplayshortskip}{6.0pt plus 3.0pt}

\setlength{\baselineskip}{13.0pt}
\setlength{\lineskip}{0.0pt}
\setlength{\lineskiplimit}{0.0pt}

%
%

\title{Mixed sheaves on Shimura varieties and their higher direct images
in toroidal compactifications \footnotemark
\footnotetext{To appear in J.\ of Alg.\ Geom.}}
\author{\footnotesize by\\ \\
\mbox{\hskip-2cm
\begin{minipage}{6cm} \begin{center} \begin{tabular}{c}
J\"org Wildeshaus\\[5pt]
\footnotesize Institut Galil\'ee\\[-3pt]
\footnotesize Universit\'e Paris 13\\[-3pt]
\footnotesize Avenue Jean-Baptiste Cl\'ement\\[-3pt]
\footnotesize F-93430 Villetaneuse\\[-3pt]
\footnotesize France\\
{\footnotesize \tt wildesh@math.univ-paris13.fr}
\end{tabular} \end{center} \end{minipage}
\hskip-2cm} 
\\[2cm] }

\date{July 28, 1999} 

\maketitle



\vfill

\noindent {\footnotesize Math.\ Subj.\ Class.\ numbers: 
14 G 35 (11 G 18, 14 D 07, 14 F 20, 
14 F 25, 19 F 27, 32 G 20).}

%

%
%
\setcounter{section}{-1}
\section{Introduction}
\label{Intro}


In this paper, we consider a {\it toroidal
compactification} of a {\it mixed Shimura variety}
\[
j: M \hookrightarrow M(\FS) \: .
\]
According to \cite{P1}, the {\it boundary} $M(\FS) - M$ has a natural
{\it stratification} into locally closed subsets,
each of which is itself (a quotient by
the action of a finite group of) a Shimura variety. Let
\[
i: M' \hookrightarrow M(\FS)
\]
be the inclusion of an individual such stratum. Both in the Hodge and the
$\ell$-adic context, there is a theory of {\it mixed sheaves}, and in particular,
a functor
\[
i^*j_*
\]
from the bounded derived category of mixed sheaves on $M$ to that of mixed sheaves on $M'$.

The objective of the present article is a formula for the effect of $i^*j_*$
on those complexes of 
mixed sheaves coming about via the {\it canonical construction},
denoted $\mu$:
The Shimura variety $M$ is associated to a linear algebraic group $P$ over $\BQ$, and any complex of algebraic representations $\BV^\bullet$ of $P$ gives rise
to a complex of
mixed sheaves $\mu(\BV^\bullet)$ on $M$. Let $P'$ be the group belonging to $M'$;
it is the quotient by a normal unipotent subgroup $U'$ of a
subgroup $P_1$ of $P$:
\[
\begin{array}{ccccc}
U' & \trianglelefteq & P_1 & \le & P  \\
& & \downarrow & & \\
& & P' & &
\end{array}
\]
Our main result (\ref{2H} in the Hodge setting; \ref{3J} in the $\ell$-adic setting) expresses the composition $i^*j_* \circ \mu$
in terms of the canonical construction $\mu'$ on $M'$, and Hochschild
cohomology of $U'$. It may be seen as complementing results 
of Harris and Zucker (\cite{HZ}), and of Pink (\cite{P2}). 

In the $\ell$-adic setting, \cite{P2} treats the analogous question for the 
natural stratification of the {\it Baily--Borel compactification} $M^*$
of a {\it pure} Shimura variety $M$. 
The resulting formula (\cite{P2}~(5.3.1)) has a more
complicated structure than ours: Besides Hochschild cohomology of a unipotent
group, it also involves cohomology of a certain arithmetic group.
Although we are interested in a different geometric situation, much of the
abstract material developed in the first two sections of \cite{P2} will enter
our proof. We should mention that the proof of Pink's result actually involves
a toroidal compactification. The stratification used is the
one induced by the stratification of $M^*$, and is therefore coarser than the
one considered in the present work.

In \cite{HZ}, Harris and Zucker study 
the {\it Hodge structure} on the boundary cohomology of
the {\it Borel--Serre compactification} of a Shimura variety. As
in \cite{P2}, toroidal compactifications enter the proof of the
main result (\cite{HZ}~(5.5.2)). It turns out to be necessary to control
the structure of $i^*j_* \circ \mu(\BV^\bullet)$ in the case
when the stratum $M'$ is minimal. 
There, the authors 
arrive at a 
description which is equivalent to ours (\cite{HZ}~(4.4.18)). 
Although they only treat the case of a pure Shimura variety, and
do not relate their result directly to representations of the group $P'$,
it is fair to say that an important part of 
the main Hodge theoretic information
entering our proof
(see (b) below) is
already contained in \cite{HZ}~(4.4). Still, our global
strategy of proof of the main comparison result \ref{2H} is different: 
We employ Saito's \emph{specialization functor}, and
a homological yoga to reduce to two seemingly weaker
comparison statements: (a) comparison for the full functor $i^*j_* \circ \mu$, but only
on the level of local systems; (b) comparison on the level of variations of
Hodge structure, but only for $\CH^0i^*j_* \circ \mu$.

%
%
It is a pleasure to thank A.~Huber, A.~Werner, D.~Blasius, C.~Deninger,
G.~Kings, 
C.~Serp\'e, J.~Steenbrink, M.~Strauch  
and T.~Wenger for useful remarks, and
G.~Weckermann for \TeX ing my manuscript. 
I am particularly grateful to R.~Pink for pointing out an error
in an earlier version of the proof of \ref{2H}.
Finally, I am indebted to the referee for
her or his helpful comments.


\myheading{Notations and Conventions:}
Throughout the whole article, we make consistent use of the language and
the main results of \cite{P1}. 

Algebraic representations of an algebraic group are finite dimensional
by definition. If a group $G$ acts on $X$, then we write $\Cent_G X$ for
the kernel of the action. If $Y$ is a subobject of $X$, then $\Stab_G Y$
denotes the subgroup of $G$ stabilizing $Y$. 

If $X$ is a variety over $\BC\,$, 
then $D^b_\con (X(\BC))$ denotes the full triangulated subcategory
of complexes of sheaves of abelian groups on $X(\BC)$ with 
constructible cohomology. The subcategory of complexes whose cohomology
is \emph{algebraically} constructible is denoted by $D^b_\con (X)$.
If $F$ is a coefficient field,
then we define triangulated categories of complexes of sheaves of
$F$-vector spaces
\[
D^b_\con (X , F) \subset D^b_\con (X(\BC) , F)
\]
in a similar fashion. The category
$\Perv_F X$ is defined as the heart of the perverse 
$t$-structure on $D^b_\con (X , F)$.

Finally, the ring of finite ad\`eles over $\BQ$ is denoted by $\BA_f$.

%
%

\section{Strata in toroidal compactifications}
\label{1}


This section provides complements to certain aspects of Pink's treatment (\cite{P1}). The first concerns the shape of the canonical stratification of a toroidal compactification of a Shimura variety. According to 
\cite{P1}~12.4~(c), these strata are quotients by finite group actions of ``smaller'' Shimura varieties. We shall show (\ref{1F}) that under mild restrictions (neatness of the compact group, and condition $(+)$ below), 
the finite groups occurring are in fact trivial.

The second result concerns the formal completion of a stratum. 
Under the above restrictions, we show (\ref{1M}) that the completion
in the toroidal compactification is canonically isomorphic to the
completion in a suitable torus embedding.
Under special assumptions on the cone decomposition giving rise to the compactification, this result is an immediate consequence of \cite{P1}~12.4~(c), which concerns the {\it closure} of the stratum in question.

Finally (\ref{1R}), we identify the normal cone of a
stratum in a toroidal compactification.\\

Let $(P, \FX)$ be {\it mixed Shimura data} (\cite{P1}~Def.~2.1). 
So in particular, $P$ is a connected algebraic linear group over $\BQ$, and
$P(\BR)$ acts on the complex manifold $\FX$ by analytic automorphisms.
Any {\it admissible parabolic subgroup} (\cite{P1}~Def.~4.5) $Q$ of $P$ has a canonical normal subgroup $P_1$ (\cite{P1}~4.7). There is a finite collection of {\it rational boundary components} $(P_1 , \FX_1)$, indexed by the $P_1 (\BR)$-orbits in $\pi_0 (\FX)$ (\cite{P1}~4.11). The $(P_1 , \FX_1)$ are themselves mixed Shimura data. 

Denote by $W$ the unipotent radical of $P$. If $P$ is reductive, i.e., if $W=0$, then $(P, \FX)$ is called {\it pure}.

Consider the following condition on $(P, \FX)$:

\begin{enumerate}
\item [$(+)$] If $G$ denotes the maximal reductive quotient of $P$, then the neutral connected component $Z (G)^0$ of the center $Z (G)$ of $G$ is, up to isogeny, a direct product of a $\BQ$-split torus with a torus $T$ of compact type (i.e., $T(\BR)$ is compact) defined over $\BQ$.
\end{enumerate}

From the proof of \cite{P1}~Cor.~4.10, one concludes:

\begin{Prop}\label{1A}
  If $(P, \FX)$ satisfies $(+)$, then so does any rational boundary component $(P_1 , \FX_1)$.
\end{Prop}

Denote by $U_1 \trianglelefteq P_1$ the ``weight $-2$'' part of $P_1$. It is abelian, normal in $Q$, and central in the unipotent radical $W_1$ of $P_1$.

Fix a connected component $\FX^0$ of $\FX$, and denote by
$(P_1 , \FX_1)$ the associated rational boundary component. There is a natural open embedding
\[
\iota: \FX^0 \longrightarrow \FX_1
\]
(\cite{P1}~4.11, Prop.~4.15~(a)). If $\FX^0_1$ denotes the connected component of $\FX_1$ containing $\FX^0$, then the image of the embedding can be described by means of the map {\it imaginary part}
\[
\imm : \FX_1 \longrightarrow U_1 (\BR) (-1) := \frac{1}{2 \pi i} \cdot U_1 (\BR) \subset U_1 (\BC)
\]
of \cite{P1}~4.14: $\FX^0$ is the preimage of an open convex cone
\[
C (\FX^0 , P_1) \subset U_1 (\BR) (-1)
\]
under $\imm |_{\FX^0_1}$ (\cite{P1}~Prop.~4.15~(b)).\\

Let us indicate the definition of the map $\imm$: given $x_1 \in \FX^0_1$, there is exactly one element $u_1 \in U_1 (\BR) (-1)$ such that $u^{-1}_1 (x_1) \in \FX^0_1$ is real, i.e., the associated morphism of the Deligne torus
\[
\iint (u^{-1}_1) \circ h_{x_1} : \BS_{\BC} \longrightarrow P_{1,\BC}
\]
(\cite{P1}~2.1) descends to $\BR$. Define $\imm (x_1) := u_1$.\\

We now describe the composition
\[
\imm \circ \iota : \FX^0 \longrightarrow U_1 (\BR) (-1)
\]
in terms of the group
\[
H_0 := \{ (z,\alpha) \in \BS \times \GL_{2,\BR} \tei N(z) = \det (\alpha) \}
\]
of \cite{P1}~4.3. Let $U_0$ denote the copy of $\BG_{a,\BR}$ in $H_0$ consisting of elements
\[
\left( 1 , \left( 
    \begin{array}{cc}
1 & \ast \\ 0 & 1
    \end{array} \right) \right) \; .
\]
According to \cite{P1}~Prop.~4.6, any $x \in \FX$ defines a morphism
\[
\omega_x : H_{0,\BC} \longrightarrow P_{\BC} \; .
\]

\begin{Lem}\label{1B}
  Let $x \in \FX^0$. Then
\[
\imm (\iota x) \in U_1 (\BR) (-1)
\]
lies in $\omega_x (U_0 (\BR) (-1) - \{0\})$.
\end{Lem}

\begin{Proof}
  Since the associations
\[
x \longmapsto \omega_x
\]
and
\[
x \longmapsto \imm (\iota x)
\]
are $(U (\BR)(-1))$-equivariant, we may assume that $\imm (x) = 0$, i.e., that
\[
h_x : \BS_{\BC} \longrightarrow P_{\BC}
\]
descends to $\BR$. According to the proof of \cite{P1}~Prop.~4.6,
\[
\omega_x : H_{0,\BC} \longrightarrow P_{\BC}
\]
then descends to $\BR$. Now
\[
h_{\iota x} : \BS_{\BC} \longrightarrow P_{1,\BC} \hookrightarrow P_{\BC}
\]
is given by $\omega_x \circ h_{\infty}$, for a certain embedding
\[
h_{\infty} : \BS_{\BC} \longrightarrow H_{0,\BC}
\]
(\cite{P1}~4.3). 

More concretely, as can be seen from \cite{P1}~4.2--4.3,
there is a $\tau \in \BC - \BR$ such that on $\BC$-valued points,
we have
\[
h_\infty: (z_1,z_2) \longto 
\left( (z_1,z_2) , \left( 
    \begin{array}{cc}
    z_1z_2 & \tau (1-z_1z_2) \\ 
    0 & 1
    \end{array} \right) \right) \; .
\]
Hence there is an element
\[
u_0 \in U_0 (\BR) (-1) - \{0\} 
\]
such that $\iint (u^{-1}_0) \circ h_\infty$ descends to $\BR$.
But then $\omega_x (u_0)$ has the defining property of $\imm (\iota x)$.
\end{Proof}

Let $F$ be a field of characteristic $0$. By definition of Shimura data, any
algebraic representation
\[
\BV \in \Rep_F P
\]
comes equipped with a natural weight filtration $W_{\bullet}$
(see \cite{P1}~Prop.~1.4). Lemma \ref{1B} 
enables us to relate it
to the weight filtration $M_{\bullet}$ of 
\[
\Res^P_{P_1} (\BV) \in \Rep_F P_1 \; :
\]

\begin{Prop}\label{1C}
  Let $\BV \in \Rep_F P$, and $T \in U_1 (\BQ)$ such that 
\[
\pm \frac{1}{2 \pi i} T \in C (\FX^0 , P_1) \; .
\] 
Then the weight filtration of $\log T$ relative to $W_{\bullet}$ (\cite{D}~(1.6.13)) exists, and
is identical to $M_{\bullet}$.
\end{Prop}

\begin{Proof}
Set $N:= \log T$. Since $\Lie (U_1)$ is of weight $-2$, we clearly have
\[
NM_i \subset M_{i-2} \; .
\]
It remains to prove that
\[
N^k : \Gr^M_{m+k} \Gr^W_m \BV \longrightarrow \Gr^M_{m-k} \Gr^W_m \BV
\]
is an isomorphism. According to \ref{1B}, there are $x \in \FX^0$ and $u_0 \in U_0 (\BR) (-1) - \{0\}$ such that
\[
\omega_x : H_{0,\BC} \longrightarrow P_{\BC}
\]
maps $u_0$ to $\pm \frac{1}{2\pi i} T$. By definition,
$M_{\bullet}$ is the weight filtration associated to the morphism
\[
\omega_x \circ h_{\infty} : \BS_{\BC} \longrightarrow P_{1,\BC} \; .
\]
Our assertion has become one about representations of $H_{0,\BC}$. But 
$\Rep_\BC H_{0,\BC}$ is semisimple, the irreducible objects being given by
\[
\Sym^n V \otimes \chi \; ,
\]
$V$ the standard representation of $\GL_{2,\BC} \; , \; \chi$ a character of $H_{0,\BC}$ and $n \ge 1$. It is straightforward to show that for any such representation, the weight filtration defined by $h_{\infty}$ equals the monodromy weight filtration for $\log u_0$.
\end{Proof}

\begin{Cor}\label{1D}
Let $T \in U_1 (\BQ)$ such that $\pm \frac{1}{2 \pi i} T \in C (\FX^0 , P_1)$.
Then
\[
\Cent_W(T) = \Cent_W(U_1) = W \cap P_1 \; .
\]
\end{Cor}

\begin{Proof}
The inclusions ``$\supset$'' hold since the right hand side is contained in $W_1$, and $U_1$ is central in $W_1$. For the reverse inclusions, let us show that
\[
\Lie \left( \Cent_W (T) \right) \subset \Lie W
\]
is contained in the (weight $\le -1$)-part of the 
restriction of the adjoint representation
\[
\Lie W  \in \Rep_\BQ P 
\]
to $P_1$. Observe that with respect to this representation, we have
\[
\ker \left( \log T \right) = \Lie \left( \Cent_W (T) \right) \; .
\]
First, recall (\cite{P1}~2.1) that $\Gr_m^{W_\bullet} (\Lie W) = 0$ 
for $m \ge 0$.
From the defining property of the weight filtration $M_{\bullet}$ of $\log T$
relative to $W_{\bullet}$, it follows that
\[
\ker \left( \log T \right) \subset M_{-1} \left( \Lie W \right) \; .
\]
Proposition~\ref{1C} guarantees that the right hand side equals the
(weight $\le -1$)-part under the action of $P_1$.
Our claim therefore follows from the equality
\[
M_{-1} \left( \Lie W \right) = \Lie \left( W \cap P_1 \right)
\]
(\cite{P1}~proof of Lemma~4.8).
\end{Proof}
 
\begin{Lem}\label{1E}
Let $P_1 \trianglelefteq Q \le P$ as before, let
$\Gamma \le Q (\BQ)$ be contained in a compact subgroup of $Q (\BA_f)$, and assume that $\Gamma$ centralizes $U_1$. Then a subgroup of finite index in $\Gamma$ is contained in
\[
(Z (P) \cdot P_1) (\BQ) \; .
\]
If $(+)$ holds for $(P , \FX)$, then a subgroup of finite index in $\Gamma$ is contained in $P_1 (\BQ)$.
\end{Lem}

\begin{Proof}
  The two statements are equivalent: if one passes from $(P , \FX)$ to the quotient data $(P, {\FX}) / Z (P)$ (\cite{P1}~Prop.~2.9), then $(+)$ holds. So assume that $(+)$ is satisfied. Fix a point $x \in \FX$, and consider the associated homomorphism
\[
\omega_x : H_{0,\BC} \longrightarrow P_{\BC} \; .
\]
Since $\omega_x$ maps the subgroup $U_0$ of $H_0$ to $U_1$, the elements in the centralizer of $U_1$ also commute with $\omega_x (U_0)$.

First assume that $P = G = G^{\ad}$. By looking at the decomposition of $\Lie G_{\BR}$ under the action of $H_0$ (\cite{P1}~Lemma~4.4~(c)), one sees that
the Lie algebra of the centralizer in $Q_{\BR}$ of $\omega_x(U_0)$,
\[
\Lie (\Cent_{Q_{\BR}} U_0) \subset \Lie Q_{\BR}
\]
is contained in $\Lie P_{1,\BR} + \Lie (\Cent_{G_{\BR}} \imm (\omega_x))$.
But $\Cent_{G_{\BR}} \imm (\omega_x)$ is a compact group, hence the image of $\Gamma$ in $(Q / P_1) (\BQ)$ is finite.

Next, if $P = G$, then by the above,
\[
\Gamma \cap (Z (G) \cdot P_1) (\BQ) 
\]
is of finite index in $\Gamma$. Because of $(+)$, the image of $\Gamma$ in $(Q / P_1) (\BQ)$ is again finite.

In the general case,
\[
\Gamma \cap (W \cdot P_1) (\BQ)
\]
is of finite index in $\Gamma$. Analysing the decomposition of $\Lie W_{\BR}$ under the action of $H_0$ (\cite{P1}~Lemma~4.4~(a) and (b)), 
or using Corollary~\ref{1D}, one realizes that
\[
\Lie (\Cent_{Q} U_1) \cap \Lie W \subset \Lie P_{1} \; .
\]
\end{Proof}

The {\it Shimura varieties} associated to mixed Shimura data $(P,\FX)$ are indexed by the open compact subgroups of $P (\BA_f)$. If $K$ is one such, then the analytic space of $\BC$-valued points of the corresponding variety $M^K := M^K (P,\FX)$ is given as
\[
M^K (\BC) := P (\BQ) \backslash ( \FX \times P (\BA_f) / K ) \; .
\]
In order to discuss compactifications, we need to introduce 
the {\it conical complex} associated to $(P,\FX)$: 
set-theoretically, it is defined as
\[
\CC (P , \FX) := \coprod_{(\FX^0 , P_1)} C (\FX^0 , P_1) \; .
\]

By \cite{P1}~4.24, the conical complex is naturally equipped with a topology (which is usually different from the coproduct topology). The closure $C^{\ast} (\FX^0 , P_1)$ of $C (\FX^0 , P_1)$ 
inside $\CC (P , \FX)$
can still be considered as a convex cone in $U_1 (\BR) (-1)$, with the induced topology. \\

For fixed $K$, the (partial) {\it toroidal compactifications} of $M^K$ are para\-me\-terized by {\it $K$-admissible partial cone decompositions}, which are collections of subsets of
\[
\CC (P,\FX) \times P (\BA_f)
\]
satisfying the axioms of \cite{P1}~6.4. If $\FS$ is one such, then in particular any member of $\FS$ is of the shape
\[
\sigma \times \{ p \} \; ,
\]
$p \in P (\BA_f)$, $\sigma \subset C^{\ast} (\FX^0 , P_1)$ a {\it convex rational polyhedral cone} in 
the vector space $U_1 (\BR) (-1)$ (\cite{P1}~5.1) not containing any non-trivial linear subspace.\\

Let $M^K (\FS) := M^K (P , \FX , \FS)$ be the associated compactification. It comes equipped with a natural stratification into locally closed strata, each of which looks as follows: Fix a pair $(\FX^0 , P_1)$ as above, $p \in P (\BA_f)$ and
\[
\sigma \times \{ p \} \in \FS
\]
such that $\sigma \subset C^{\ast} (\FX^0 , P_1)$. Assume that
\[
\sigma \cap C (\FX^0 , P_1) \neq \emptyset \; .
\]
To $\sigma$, one associates Shimura data
\[
\left( \Pes , \Xes \right)
\]
(\cite{P1}~7.1), whose underlying group $\Pes$ is the quotient of $P_1$ by the algebraic subgroup
\[
\langle \sigma \rangle \subset U_1
\]
satisfying $\BR \cdot \sigma = \frac{1}{2 \pi i} \cdot \langle \sigma \rangle (\BR)$. Set 
\[
K_1  :=  P_1 (\BA_f) \cap p \cdot K \cdot p^{-1} \; , \quad
\pi_{[\sigma]}  :  P_1 \longonto \Pes \; .
\]
According to \cite{P1}~7.3, there is a canonical map
\[
i (\BC) : \Mps (\Pes , \Xes) (\BC) \longrightarrow 
M^K (\FS) (\BC)
:= M^K (P , \FX , \FS) (\BC)
\]
whose image is locally closed. In fact, $i (\BC)$ is a quotient map onto its image.

\begin{Prop}\label{1F}
  Assume that $(P, \FX)$ satisfies $(+)$, and that $K$ is \emph{neat} (see e.g.\ \cite{P1}~0.6). Then $i (\BC)$ is injective, i.e., it identifies $\Mps (\BC)$ with a locally closed subspace of $M^K (\FS) (\BC)$.
\end{Prop}

\begin{Proof}
Consider the group $\Delta_1$ of \cite{P1}~6.18:
 \begin{eqnarray*}
   H_Q & := & \Stab_{Q (\BQ)} (\FX_1) \cap P_1 (\BA_f) \cdot p \cdot K \cdot p^{-1} \; , \\
\Delta_1 & := & H_Q / P_1 (\BQ) \; .
 \end{eqnarray*}
The subgroup $\Delta_1 \le (Q / P_1) (\BQ)$ is arithmetic. According to \cite{P1}~7.3, the image under $i (\BC)$ is given by the quotient of $\Mps (\BC)$ by a certain subgroup
\[
\Stab_{\Delta_1} ([\sigma]) = \Stab_{H_Q} ([\sigma]) / P_1 (\BQ) \le \Delta_1 \; .
\]
This stabilizer refers to the action of $H_Q$ on the double quotient
\[
P_1 (\BQ) \backslash \FS_1 / P_1 (\BA_f)
\]
of \cite{P1}~7.3. Denote the projection $Q \to Q / P_1$ by $\pr$, so $\Delta_1 = \pr (H_Q)$, and
\[
\Stab_{\Delta_1} ([\sigma]) = \pr  \left( \Stab_{H_Q} ([\sigma]) \right) \; .
\]
By Lemma~\ref{1G}, this group is trivial under the hypotheses of
the proposition.
\end{Proof}

\begin{Lem}\label{1G}
  If $(P , \FX)$ satisfies $(+)$ then $\Stab_{\Delta_1} ([\sigma])$ is finite. If, in addition, $K$ is neat then $\Stab_{\Delta_1} ([\sigma]) = 1$.
\end{Lem}

\begin{Proof}
  The second claim follows from the first since $\Stab_{\Delta_1} ([\sigma])$ is contained in
\[
(Q / P_1) (\BQ) \cap \pr (p \cdot K \cdot p^{-1}) \; ,
\]
which is neat if $K$ is.

Consider the arithmetic subgroup of $Q (\BQ)$
\[
\Gamma_Q := H_Q \cap p \cdot K \cdot p^{-1} \; .
\]
The group $\pr (\Gamma_Q)$ is arithmetic, hence of finite index in $\Delta_1$. Hence
\[
\Stab_{\pr(\Gamma_Q)}([\sigma]) = 
\pr \left( \Stab_{\Gamma_Q} ([\sigma])\right) \le \Stab_{\Delta_1} ([\sigma])
\]
is of finite index. Now
\[
\Stab_{\Gamma_Q} (\sigma) \le \Stab_{\Gamma_Q} ([\sigma])
\]
is of finite index. By \cite{P1}~Thm.~6.19, a subgroup of finite index of $\Stab_{\Gamma_Q} (\sigma)$ centralizes $U_1$. Our claim thus follows from 
Lemma~\ref{1E}.
\end{Proof}

\begin{Rem}\label{1H}
  The lemma shows that the groups ``\, $\Stab_{\Delta_1} ([\sigma])$'' occurring in 7.11, 7.15, 7.17, 9.36, 9.37, and 12.4 of \cite{P1} are all trivial provided that $(P,\FX)$ satisfies $(+)$ and $K$ is neat.
\end{Rem}

We continue the study of the map
\[
i (\BC) : \Mps (\BC) \longrightarrow M^K (\FS) (\BC) \; .
\]
Let $\Ses$ be the minimal $K_1$-admissible cone decomposition of
\[
\CC (P_1 , \FX_1) \times P_1 (\BA_f)
\]
containing $\sigma \times \{ 1 \}$; $\Ses$ can be realized inside the decomposition $\FS^0_1$ of \cite{P1}~6.13; by definition, it is {\it concentrated in the unipotent fibre} (\cite{P1}~6.5~(d)). 
View $\Mps (\BC)$ as sitting inside $M^{K_1} (\Ses) (\BC)$:
\[
i_1 (\BC) : \Mps (\BC) \hookrightarrow M^{K_1} (\Ses) (\BC) \; .
\]
Consider the diagram
\[
\vcenter{\xymatrix@R-10pt{ 
\Mps (\BC) \ar@{^{ (}->}[r]^{i_1 (\BC)} \ar@{_{ (}->}[dr]_{i(\BC)} &
M^{K_1} (\Ses) (\BC) \\
&  M^K (\FS) (\BC) \\}}
\]
\cite{P1}~6.13 contains the definition of an open neighbourhood
\[
\FU := \OFU (P_1 , \FX_1 , p)
\]
of $\Mps (\BC)$ in $M^{K_1} (\Ses) (\BC)$, and a natural extension $f$ of the map $i (\BC)$ to $\FU\,$:
\[
\vcenter{\xymatrix@R-10pt{ 
\Mps (\BC) \ar@{^{ (}->}[r] \ar@{_{ (}->}[dr]_{i (\BC)} & \FU \ar[d]^f \\
&  M^K (\FS) (\BC) \\}}
\]
\begin{Prop}\label{1I}
(a) $f$ is open.\\
(b) We have the equality
\[
f^{-1} (M^K(\BC)) = \FU \cap M^{K_1} (\BC) \; .
\]
\end{Prop}

\begin{Proof} Let us recall the definition of $\OFU (P_1 , \FX_1,p)$, and part of the construction of 
$M^K (\FS) (\BC)$:
Let $\FX^+ \subset \FX$ be the preimage under
\[
\FX \longrightarrow \pi_0 (\FX)
\]
of the $P_1 (\BR)$-orbit in $\pi_0 (\FX)$ associated to $\FX_1$, and
\[
\FX^+ \longrightarrow \FX_1
\]
the map discussed after Proposition~\ref{1A}; according to \cite{P1}~Prop.~4.15~(a), it is still an open embedding (i.e., injective). 
As in \cite{P1}~6.10, set
\[
\FU (P_1,  \FX_1 , p) := P_1 (\BQ) \backslash 
( \FX^+ \times P_1 (\BA_f) / K_1) \; .
\]
It obviously admits an open embedding into $M^{K_1} (\BC)$ as well as an open morphism to $M^K (\BC)$. One defines (\cite{P1}~6.13)
\[
\OFU (P_1 , \FX_1 , p) \subset M^{K_1} (\Ses) (\BC)
\]
as the interior of the closure of $\FU (P_1 , \FX_1 , p)$. 
Then $M^K (\FS) (\BC)$ is defined as the quotient with respect to some equivalence relation $\sim$ on the disjoint sum of all $\OFU (P_1 , \FX_1 , p)$ (\cite{P1}~6.24). In particular, for our {\it fixed} choice of $(P_1 , \FX_1)$ and $p$, there is a continuous map
\[
f : \OFU (P_1 , \FX_1 , p) \longrightarrow M^K (\FS) (\BC) \; .
\]
From the description of $\sim$ (\cite{P1}~6.15--6.16), one sees that $f$ is open; the central point is that the maps
\[
\overline{\beta} := \overline{\beta} (P_1 , \FX_1 , P'_1 , \FX'_1 , p):
\OFU (P_1,  \FX_1 , p) \cap M^{K_1} (P_1 , \FX_1, \FS''^0) (\BC) \longto
\OFU (P'_1,  \FX'_1 , p) 
\]
of \cite{P1}~page~152 are open. This shows (a). As for (b), one observes that
\[
\overline{\beta}^{-1} \left( \FU (P'_1,  \FX'_1 , p) \right) =
\FU (P_1,  \FX_1 , p) \; .
\]
\end{Proof}

\begin{Rem}\label{1J} 
\cite{P1}~Cor.~7.17 gives a much stronger statement than Proposition~\ref{1I}~(a), assuming that $\FS$ is \emph{complete} (\cite{P1}~6.4) and satisfies condition ($\ast$) of \cite{P1}~7.12. In this case, one can identify a suitable open neighbourhood of the \emph{closure} of
\[
\Stab_{\Delta_1} ([\sigma]) \backslash \Mps (\BC) = \image( i(\BC) ) \subset M^K (\FS) (\BC)
\]
with an open neighbourhood of the closure of
\[
\Stab_{\Delta_1} ([\sigma]) \backslash \Mps (\BC) \subset \Stab_{\Delta_1} ([\sigma]) \backslash M^{K_1} (\FS_1) (\BC) \; ,
\]
where 
\[
\Ses \subset \FS_1 := ([\cdot p]^{\ast} \FS) \, |_{(P_1 , \FX_1)} 
\]
(\cite{P1}~6.5~(a) and (c)).

Consequently, one can identify the formal completions (in the sense of analytic spaces) of $M^K (\FS) (\BC)$ and of
\[
\Stab_{\Delta_1} ([\sigma]) \backslash M^{K_1} (\FS_1) (\BC)
\]
along the closure of the stratum
\[
\Stab_{\Delta_1} ([\sigma]) \backslash \Mps (\BC) \; .
\]
\end{Rem}

It will be important to know that without the hypotheses of
\cite{P1}~Cor.~7.17, the completions along the stratum itself still agree. For simplicity, we assume that the hypotheses of Proposition~\ref{1F} are met, and hence that $\Stab_{\Delta_1} ([\sigma]) = 1$.

\begin{Thm} \label{1K}
Assume that $(P, \FX)$ satisfies $(+)$, and that $K$ is neat. 
\begin{itemize}
\item[(i)] The map $f$ of \ref{1I} is locally biholomorphic near $\Mps (\BC)$.
\item[(ii)] 
$f$ induces an isomorphism between the formal analytic completions of
$M^K (\FS) (\BC)$ and of $M^{K_1} (\Ses) (\BC)$ along $\Mps (\BC)$. 
\end{itemize}
\end{Thm}

\begin{Proof}
  $f$ is open and identifies the analytic subsets
\[
\Mps (\BC) \subset M^{K_1} (\Ses) (\BC)
\]
and
\[
\Mps (\BC) \subset M^K (\FS) (\BC) \; .
\]
For (ii), we have to compare certain sheaves of functions. The claim therefore follows from (i).

According to \cite{P1}~6.18, the image of $f$ equals the quotient of $\FU$ by the action of a group $\Delta_1$ of analytic automorphisms, which according to \cite{P1}~Prop.~6.20 is properly discontinuous.
\end{Proof}

So far, we have worked in the category of analytic spaces. According to 
Pink's generalization to mixed Shimura varieties of the Algebraization
Theorem of Baily and Borel
(\cite{P1}~Prop.~9.24), there exist canonical structures of normal algebraic varieties on the $M^K (P , \FX) (\BC)$, which we denote as
\[
M^K_{\BC} := M^K (P, \FX)_{\BC} \; .
\]
If there exists a structure of normal algebraic variety on $M^K (P , \FX , \FS) (\BC)$ extending $M^K_{\BC}$, then it is unique (\cite{P1}~9.25); denote it as
\[
M^K (\FS)_{\BC} := M^K (P , \FX , \FS)_{\BC} \; .
\]
Pink gives criteria on the existence of $M^K (\FS)_{\BC}$ (\cite{P1}~9.27, 9.28). If any cone of a cone decomposition $\FS'$ for $(P , \FX)$ is contained in a cone of $\FS$, and both $M^K (\FS')_{\BC}$ and $M^K (\FS)_{\BC}$ exist, then the morphism
\[
M^K (\FS') (\BC) \longrightarrow M^K (\FS) (\BC)
\]
is algebraic (\cite{P1}~9.25). From now on we implicitly assume the existence whenever we talk about $M^K (\FS)_{\BC}$.

According to \cite{P1}~Prop.~9.36, the stratification of $M^K (\FS)_{\BC}$ holds algebraically.

\begin{Thm}\label{1L}
  Assume that $(P, \FX)$ satisfies $(+)$, and that $K$ is neat. The isomorphism of Theorem~\ref{1K} induces a canonical isomorphism between the 
formal completions of $M^K (\FS)_{\BC}$ and of $M^{K_1} (\Ses)_{\BC}$ along $\MpC$. 
\end{Thm}

\begin{Proof}  
If $\FS$ is complete and satisfies ($\ast$) of \cite{P1}~7.12, then this is an immediate consequence of \cite{P1}~Prop.~9.37, which concerns the formal completions along the closure of $\MpC$.

We may replace $K$ by a normal subgroup $K'$ of finite index: the objects on
the level of $K$ come about as quotients under the finite group $K / K'$
of those on the level of $K'$. Therefore,
we may assume, thanks to \cite{P1}~Prop.~9.29 and Prop.~7.13, that there is a complete cone decomposition $\FS'$ containing $\sigma \times \{ p \}$ and satisfying ($\ast$) of \cite{P1}~7.12. 
Let $\FS''$ be the coarsest refinement of both $\FS$ and $\FS'$; it still contains $\sigma \times \{ p \}$, and $M^K (\FS'')_{\BC}$ exists because of \cite{P1}~Prop.~9.28. We have
\[
\Ses = \FS''_{1,[\sigma]} = \FS'_{1,[\sigma]} \; ,
\]
hence the formal completions all agree analytically. But on the level of $\FS'_{1,[\sigma]}$, the isomorphism is algebraic.
\end{Proof}

According to \cite{P1}~Thm.~11.18, there exists a {\it canonical model} of the variety $M^K (P , \FX)_{\BC}$, which we denote as
\[
M^K := M^K (P, \FX) \; .
\]
It is defined over the {\it reflex field} $E (P , \FX)$ of $(P , \FX)$ (\cite{P1}~11.1). The reflex field does not change when passing from $(P , \FX)$ to a rational boundary component (\cite{P1}~Prop.~12.1).

If $M^K (\FS)_{\BC}$ exists, then it has a canonical model $M^K (\FS)$ over $E (P , \FX)$ extending $M^K$, and the stratification descends to $E (P, \FX)$. In fact, \cite{P1}~Thm.~12.4 contains these statements under special hypotheses on $\FS$. However, one passes from $\FS$ to a covering by finite cone decompositions (corresponding to an open covering of $M^K (\FS)_{\BC}$), and then (\cite{P1}~Cor.~9.33) to a subgroup of $K$ of finite index to see that the above claims hold as soon as $M^K (\FS)_{\BC}$ exists.

\begin{Thm}\label{1M}
  Assume that $(P , \FX)$ satisfies $(+)$, and that $K$ is neat. The isomorphism of Theorem~\ref{1L} descends to a canonical isomorphism between the formal completions of $M^K (\FS)$ and of $M^{K_1} (\Ses)$ along $\Mps$.
\end{Thm}

\begin{Proof}
  If $\FS$ is complete and satisfies ($\ast$) of \cite{P1}~7.12, then this statement is contained in \cite{P1}~Thm.~12.4~(c).

In fact, the proof of \cite{P1}~Thm.~12.4~(c) does not directly use the special hypotheses on $\FS$: the strategy is really to prove \ref{1M} and then deduce the stronger conclusion of \cite{P1}~12.4~(c) from the fact that it holds over $\BC$\,; the point there is (\cite{P1}~12.6) that since the schemes are normal, morphisms descend if they descend on some open dense subscheme.

Thus the proof of our claim is contained in \cite{P1}~12.7--12.17.
\end{Proof}

\begin{Rem}\label{1O} (a) Without any hypotheses on $(P , \FX)$ and $K$, there are obvious variants of Theorems~\ref{1K}, \ref{1L}, and \ref{1M}. In particular, there is a canonical isomorphism between the formal completions of $M^K (\FS)$ and of
\[
\Stab_{\Delta_1} ([\sigma]) \backslash M^{K_1} (\Ses)
\]
along
\[
\Stab_{\Delta_1} ([\sigma]) \backslash \Mps \; .
\]
(b) By choosing simultaneous refinements, one sees that the isomorphisms of \ref{1K}~(ii), \ref{1L}, and \ref{1M} do not depend on the cone decomposition $\FS$ ``surrounding'' our fixed cone $\sigma \times \{ p \}$.
\end{Rem}

In the situation
we have been considering, the cone $\sigma$ is called {\it smooth} 
with respect to the lattice
\[
\Gamma^p_U (-1) := \frac{1}{2 \pi i} \cdot \left( U_1 (\BQ) \cap K_1 \right) \subset \frac{1}{2 \pi i} \cdot U_1 (\BR) 
\]
if the semi-group
\[
\Lambda_{\sigma} := \sigma \cap \Gamma^p_U (-1)
\]
can be generated (as semi-group) 
by a subset of a $\BZ$-basis of $\Gamma^p_U (-1)$. The corresponding statement is then necessarily true for any face of $\sigma$. Hence the $K_1$-admissible partial cone decomposition $\Ses$ is smooth in the sense of \cite{P1}~6.4.\\

Let us introduce the following condition on $(P_1,\FX_1)$:

\begin{enumerate}
\item [$(\cong)$] The canonical morphism
$(\Pes,\Xes) \longto (P_1,\FX_1) / \langle \sigma \rangle$
(\cite{P1}~7.1) is an isomorphism.
\end{enumerate}

In particular, there is an epimorphism of Shimura data from $(P_1,\FX_1)$ to $(\Pes,\Xes)$. According to \cite{P1}~7.1, we have:

\begin{Prop} \label{1P}
Condition $(\cong)$ is satisfied whenever $(P_1,\FX_1)$ is a \emph{proper}
boundary component of some other mixed Shimura data, e.g., if the
parabolic subgroup $Q \le P$ is proper.
\end{Prop}

Under the hypothesis $(\cong)$, we can establish more structural
properties of our varieties:

\begin{Lem} \label{1Q}
Assume that $(\cong)$ is satisfied.
\begin{itemize}
\item [(i)] The Shimura variety $M^{K_1}$ is a torus torsor over
$\Mps$:
\[
\pi_{[\sigma]} : M^{K_1} \longto \Mps \; .
\]
The compactification $M^{K_1} (\Ses)$ is a \emph{torus embedding} along
the fibres of $\pi_{[\sigma]}$:
\[
\overline{\pi_{[\sigma]}} : M^{K_1} (\Ses) \longto \Mps 
\]
admitting only one closed stratum. The section
\[
i_1 : \Mps \hookrightarrow M^{K_1} (\Ses)
\]
of $\overline{\pi_{[\sigma]}}$ identifies the base with this closed stratum.
\item [(ii)] Assume that $\sigma$ is smooth. Then 
\[
\overline{\pi_{[\sigma]}} : M^{K_1} (\Ses) \longto \Mps 
\]
carries a canonical structure of vector bundle, with zero section $i_1$.
The rank of this vector bundle is equal to the dimension of $\sigma$.
\end{itemize}
\end{Lem}

\begin{Proof}
(i) This is \cite{P1}~remark on the bottom of page~165, taking into account 
that $\Ses$ is minimal with respect to the property of
containing $\sigma$.\\
(ii) If $\sigma$ is smooth of dimension $c$, then by definition, the semi-group
$\Lambda_{\sigma}$ can be generated by an appropriate basis of the ambient 
real vector space. One shows that each choice of such a basis gives rise to
the same $\Mps$-linear structure on $M^{K_1} (\Ses)$.
\end{Proof}

We conclude the section by putting together all the results obtained so far:

\begin{Thm} \label{1R}
Assume that $(P,\FX)$ satisfies $(+)$, that $(P_1,\FX_1)$ satisfies $(\cong)$,
that $K$ is neat, and that $\sigma$ is smooth. Then there is a canonical isomorphism of vector bundles over $\Mps$
\[
\iota_\sigma: M^{K_1}(\Ses) \isoto N_{\Mps / M^K (\FS)}
\]
identifying $M^{K_1}(\Ses)$ and the normal bundle of $\Mps$ in $M^K (\FS)$.
\end{Thm}

\begin{Proof}
The isomorphism of Theorem~\ref{1M} induces an isomorphism
\[
N_{\Mps / M^{K_1}(\Ses)} \isoto N_{\Mps / M^K (\FS)} \; .
\]
But the normal bundle of the zero section in a vector bundle is canonically
isomorphic to the vector bundle itself.
\end{Proof}

%
%

\section{Higher direct images for Hodge modules}
\label{2}


Let $M^K(\FS) = M^K(P,\FX, \FS)$ be a toroidal compactification
of a Shimura variety $M^K = M^K(P,\FX)$, and $\Mps = \Mps(\Pes,\Xes)$
a boundary stratum. Consider the situation
\[
M^K \stackrel{j}{\hookrightarrow} M^K (\FS) 
\stackrel{i}{\hookleftarrow} \Mps  \; .
\]

Saito's formalism (\cite{Sa}) gives a functor $i^{\ast} j_{\ast}$ between the bounded derived categories of {\it algebraic mixed Hodge modules} on $M^K_{\BC}$ and on $\MpC$ respectively. The main result of this section 
(Theorem~\ref{2H}) gives a formula for the restriction of $i^{\ast} j_{\ast}$ onto the image of the natural functor associating to an algebraic representation of $P$ a variation of Hodge structure on $M^K_{\BC}$. The proof has two steps: first, one employs the \emph{specialization functor} \`a la Verdier--Saito,
and Theorem~\ref{1R}, to reduce from the toroidal to a toric situation
(\ref{2K}). The second step consists in proving the compatibility statement on the level of $\CH^0$ and then appealing to homological algebra, which implies 
compatibility on the level of functors between derived categories. \\

Throughout the whole section, we fix a set of data satisfying
the hypotheses of Theorem~\ref{1R}. We thus have
Shimura data $(P,\FX)$ satisfying condition $(+)$, a rational boundary component $(P_1 , \FX_1)$ satisfying condition $(\cong)$,
an open, compact and neat subgroup $K \le P (\BA_f)$, an element 
$p \in P (\BA_f)$ and a smooth 
cone $\sigma \times \{ p \} \subset C^{\ast} (\FX^0 , P_1) \times \{ p \}$ belonging to some $K$-admissible partial cone decomposition $\FS$ such that $M^K (\FS)$ exists. We assume that
\[
\sigma \cap C (\FX^0 , P_1) \neq \emptyset \; ,
\]
and write $K_1 := P_1 (\BA_f) \cap p \cdot K \cdot p^{-1}$,
\[
j : M^K \hookrightarrow M^K (\FS) \; ,
\]
and 
\[
i : \Mps \hookrightarrow M^K (\FS) \; .
\]
Similarly, write
\[
j_1 : M^{K_1} \hookrightarrow M^{K_1} (\Ses) \; ,
\]
and
\[
i_1 : \Mps \hookrightarrow M^{K_1} (\Ses)
\]
for the immersions into the torus embedding $M^{K_1} (\Ses)$, which according to
Theorem~\ref{1R} we identify with the normal bundle of $\Mps$ in $M^K (\FS)$.

If we denote by $c$ the dimension of $\sigma$, then both $i$ and $i_1$ are of
pure codimension $c$.\\

The immersions $i (\BC)$ and $i_1 (\BC)$ factor as 
\[
\vcenter{\xymatrix@R-10pt{ 
\Mps (\BC) \ar@{^{ (}->}[r] \ar@{=}[d] &
\FU \ar@{^{ (}->}[r] \ar[d]^f &
M^{K_1} (\Ses) (\BC) \ar@{<-^{ )}}[r] &
M^{K_1} (\BC) \\
\Mps (\BC) \ar@{^{ (}->}[r] &
\FV \ar@{^{ (}->}[r] &
M^K (\FS) (\BC) \ar@{<-^{ )}}[r] &
M^K(\BC) \\}}
\]
where $\FU$ and $\FV$ are open subsets of the respective compactifications,
and $f$ is the map of \ref{1I}. For a sheaf $\CF$ on $M^K (\BC)$, we can consider the restriction $f^{-1} \CF$ on 
$f^{-1} (M^K (\BC)) = \FU \cap M^{K_1}(\BC)$.

Let $F$ be a coefficient field
of characteristic $0$. 
Denote by 
\[
\mu_{K,\topp} : \Rep_F P \longto \Loc_F  M^K (\BC)
\]
the exact tensor functor associating to an algebraic representation $\BV$
the sheaf of
sections of
\[
P (\BQ) \backslash \left( \FX \times \BV \times P (\BA_f) / K \right) 
\]
on
\[
M^K (\BC) = P(\BQ) \backslash \left( \FX \times P (\BA_f) / K \right)\; .
\]

\begin{Prop}\label{2C}
Let $\BV \in \Rep_F P$. Then $f^{-1} \circ \mu_{K,\topp} \BV$ 
is the restriction to $f^{-1} (M^K (\BC))$ of the local system 
$\mu_{K_1,\topp} \Res^P_{P_1} \BV$ on
\[
M^{K_1} (\BC) = P_1 (\BQ) \backslash (\FX_1 \times P_1 (\BA_f) / K_1) \; .
\]
\end{Prop}

\begin{Proof}
$f^{-1} (M^K (\BC))$ equals the set
\[
\FU (P_1 , \FX_1 , p) := P_1 (\BQ) \backslash \left( \FX^+ \times P_1 (\BA_f) / K_1\right) \subset M^{K_1} (\BC)\; ,
\]
and $f \, |_{\FU (P_1 , \FX_1 , p)}$ is given by
\[
[(x , p_1)] \longmapsto [(x , p_1 \cdot p)]
\]
(\cite{P1}~6.10).
\end{Proof}

Using Theorem~\ref{1K}~(i), it is easy to construct a canonical isomorphism
of functors
\[
i^{\ast}  j_{\ast} \quad , \quad i^{\ast}_1  (j_1)_{\ast} \circ f^{-1}:
D^b_\con ( M^K (\BC) ) \longrightarrow D^b_\con ( \Mps (\BC) ) \; . 
\]
For us, it will be necessary to establish a connection between
$j_{\ast}$ and $(j_1)_{\ast} \circ f^{-1}$. This
relation will be given by Verdier's specialization functor
(\cite{V2}~9)
\[
\Ss := \SP_{\Mps} : D^b_\con ( M^K (\FS) (\BC) ,F ) \longrightarrow
                    D^b_\con ( M^{K_1}(\Ses) (\BC) , F ) \; .
\]
According to \cite{V2}~p.~358, the functor $\Ss$ has the properties
(SP0)--(SP6) of \cite{V2}~8, \emph{convenablement transpos\'ees}.
In particular:
\begin{itemize}
\item[(SP0)] It can be computed locally.
\item[(SP1)] The complexes in the image of $\Ss$ are \emph{monodromic}, i.e.,
their cohomology objects are locally constant on each $\BC^*$-orbit in
$M^{K_1}(\Ses)(\BC)$.
\item[(SP5)] We have the equality $i^* = i_1^* \circ \Ss$.
\end{itemize}
From Theorem~\ref{1K}~(i) and from (SP0), we conclude 
that in order to compute the effect of $\Ss$ on a complex 
of sheaves $\CF^\bullet$, we may pass to the complex $f^{-1} \CF^\bullet$.

On the other hand, in the case when $P=P_1$, one considers the specialization
functor for the zero section in a vector bundle. Using the definition of
$\Ss$, and hence, of the nearby cycle
functor $\psi_\pi$ in the analytic context
(\cite{Comp}~1.2), one sees that in this case, the functor $\Ss$
induces the identity on the
category of monodromic complexes.\\

By extension by zero, let us view objects of $\Loc_F M^K (\BC)$ as 
sheaves on $M^K (\FS) (\BC)$. From the above, one concludes that the functor
$\Ss$ induces a functor
\[
\Loc_F M^K (\BC) \longto \Loc_F M^{K_1} (\BC) \; ,
\]
equally denoted by $\Ss$. For local systems in the image of
$\mu_{K,\topp}$, we have:

\begin{Prop} \label{2Ca}
There is a commutative diagram
\[
\vcenter{\xymatrix@R-10pt{ 
\Rep_F P \ar[r]^{\Res^P_{P_1}} \ar[d]_{\mu_{K,\topp}} &
\Rep_F P_1                     \ar[d]^{\mu_{K_1,\topp}} \\
\Loc_F M^K (\BC) \ar[r]^{\Ss}  &
\Loc_F M^{K_1} (\BC)  \\}}
\]
\end{Prop} 

\begin{Prop} \label{2Cb}
There is a commutative diagram
\[
\vcenter{\xymatrix@R-10pt{ 
\Rep_F P \ar[r]^{\Res^P_{P_1}} \ar[d]_{\mu_{K,\topp}} &
\Rep_F P_1                     \ar[d]^{\mu_{K_1,\topp}} \\
\Loc_F M^K (\BC) \ar[d]_{j_*}  &
\Loc_F M^{K_1} (\BC) \ar[d]^{(j_1)_*}  \\
D^b_\con \left( M^K (\FS) (\BC) , F \right) \ar[r]^{\Ss} &
D^b_\con \left( M^{K_1}(\Ses) (\BC) , F  \right) \\}}
\]
\end{Prop}

Consequently:

\begin{Thm}\label{2D}
  There is a commutative diagram
\[
\vcenter{\xymatrix@R-10pt{ 
D^b \left( \Rep_F P \right) \ar[r]^{\Res^P_{P_1}} \ar[d]_{\mu_{K,\topp}} &
D^b \left( \Rep_F P_1 \right) \ar[r]^{R (\;)^{\langle \sigma \rangle}} &
D^b \left( \Rep_F P_{1,[\sigma]} \right) 
                        \ar[d]^{\mu_{\pi_{[\sigma]} (K_1),\topp}} \\
D^b_\con \left( M^K (\BC) ,F \right) \ar[rr]^{i^{\ast} j_{\ast}} &&
D^b_\con \left( \Mps (\BC) , F \right) \\}}
\]
Here, $R (\;)^{\langle \sigma \rangle}$ refers to Hochschild cohomology of the unipotent group $\langle \sigma \rangle \le P_1$.
\end{Thm}

\begin{Proof}
By (SP5) and Proposition~\ref{2Cb}, we may assume $P=P_1$.
Denote by $L_{\sigma}$ the monodromy group of $\Mps (\BC)$ inside $M^{K_1} (\Ses) (\BC)$.
It is generated by the semi-group
\[
\Lambda_{\sigma} (1) := 2 \pi i \cdot \Lambda_{\sigma} \subset U_1 (\BQ) 
\]
(see the definition before \ref{1P}),
and forms a lattice inside $\langle \sigma \rangle$.
It is well known that on the image of $\mu_{K,\topp}$, the functor 
$(i_1)^* (j_1)_*$
can be computed via group cohomology of the abstract group $L_{\sigma}$.
Since $\langle \sigma \rangle$ is unipotent, its Hochschild cohomology
coincides with cohomology of $L_\sigma$ on algebraic representations.
\end{Proof}

Let us reformulate Theorem~\ref{2D} in the language of perverse sheaves
(\cite{BBD}). Since local systems on the space of $\BC$-valued points
of a smooth complex variety can be viewed
as perverse sheaves (up to a shift), we may consider $\Mut$ as exact functor
\[
\Rep_F P \longto \Pk \; .
\]
By \cite{B}~Main Theorem~1.3, the bounded derived category 
\[
D^b \left( \Pk \right)
\]
is canonically isomorphic to $D^b_\con (M^K_\BC , F)$. Theorem \ref{2D} acquires
the following form:

\begin{Var}\label{2F} 
There is a commutative diagram
\[
\vcenter{\xymatrix@R-10pt{ 
D^b \left( \Rep_F P \right) \ar[r]^{\Res^P_{P_1}} \ar[d]_{\mu_{K,\topp}} &
D^b \left( \Rep_F P_1 \right) \ar[r]^{R (\;)^{\langle \sigma \rangle}} &
D^b \left( \Rep_F P_{1,[\sigma]} \right) 
                        \ar[d]^{\mu_{\pi_{[\sigma]} (K_1),\topp}} \\
D^b \left( \Pk \right) \ar[rr]^{i^{\ast} j_{\ast}[-c]} &&
D^b \left( \PK \right) \\}}
\]
\end{Var}

By definition of Shimura data, there is a tensor functor associating to an algebraic $F$-representation $\BV$ of $P$, for $F \subseteq \BR$, a variation of Hodge structure $\mu (\BV)$ on $\FX$ (\cite{P1}~1.18). 
It descends to a variation $\mu_K (\BV)$ on $M^K (\BC)$ with underlying local system $\mu_{K,\topp} (\BV)$.
We refer to the functor $\mu_K$ as the {\it canonical construction} of variations of Hodge structure from representations of $P$.

By \cite{W}~Thm.~2.2, the image of $\mu_K$ is contained in the category 
$\VMHS_F M^K_{\BC}$ of {\it admissible} variations, and hence
(\cite{Sa}~Thm.~3.27), in the category $\MHM_F M^K_{\BC}$ of algebraic mixed Hodge modules.\\

According to \cite{Sa}~2.30, there
is a Hodge theoretic variant of the specialization functor:
\[
\Ss := \SP_{\Mps} : \MHM_F M^K (\FS)_\BC \longto \MHM_F M^{K_1} (\Ses)_\BC \; ,
\]
which is compatible with Verdier's functor discussed earlier. Since 
the latter maps local
systems on $M^K (\BC)$ to local systems on $M^{K_1} (\BC)$ (viewed as 
sheaves on the respective compactifications by extension by zero), we see that
$\Ss$ induces a functor
\[
\VMHS_F M^K_{\BC} \longto \VMHS_F M^{K_1}_\BC \; ,
\]
equally denoted by $\Ss$. 

\begin{Thm} \label{2K}
There is a commutative diagram
\[
\vcenter{\xymatrix@R-10pt{ 
\Rep_F P \ar[r]^{\Res^P_{P_1}} \ar[d]_{\mu_K} &
\Rep_F P_1                     \ar[d]^{\mu_{K_1}} \\
\VMHS_F M^K_\BC \ar[r]^{\Ss} &
\VMHS_F M^{K_1}_\BC \\}}
\]
which is compatible with that of \ref{2Ca}.
\end{Thm}

\begin{Proof}
For $\BV \in \Rep_F P$, denote by $\BV_P$ and $\BV_{P_1}$ the two variations on the open subset $f^{-1} (M^K (\BC))$ of $M^{K_1} (\BC)$ obtained by restricting $\mu_K (\BV)$ and $\mu_{K_1} \left( \Res^P_{P_1} (\BV) \right)$ respectively.
By Proposition~\ref{2C}, the underlying local systems are identical. 

By \cite{P1}~Prop.~4.12,
the Hodge filtrations of $\BV_P$ and $\BV_{P_1}$ coincide. 

Denote the weight filtration on the variation
$\BV_P$ by $W_{\bullet}$, and that on $\BV_{P_1}$ by $M_{\bullet}$. Denote by $L_{\sigma} \subset U_1 (\BQ)$ the monodromy group of $\Mps (\BC)$ inside $M^{K_1} (\Ses) (\BC)$.
Let $T \in L_{\sigma}$ such that $\frac{1}{2 \pi i} T$ or $- \frac{1}{2 \pi i} T$ lies in $C (\FX^0 , P_1)$. According to Proposition~\ref{1C}, 
the weight filtration of $\log T$ relative to $W_{\bullet}$
is identical to $M_{\bullet}$.

Choosing $T$ as the product of the generators of the semi-group
\[
\Lambda_{\sigma}(1) \subset L_{\sigma} \; ,
\]
one concludes that $\BV_{P_1}$ carries
the limit Hodge structure of $\BV_P$ near $\Mps$. 
Using the definition of $\Ss$, and hence, of the nearby cycle functor
in the Hodge theoretic context (\cite{Sa}~2.3), one sees that the two
variations $\mu_{K_1} \circ \Res^P_{P_1} \BV$ and $\Ss \circ \mu_K \BV$ 
coincide.
\end{Proof}

\begin{Cor} \label{2Ka}
There is a commutative diagram
\[
\vcenter{\xymatrix@R-10pt{ 
\Rep_F P \ar[r]^{\Res^P_{P_1}} \ar[d]_{\mu_K} &
\Rep_F P_1                     \ar[d]^{\mu_{K_1}} \\
\VMHS_F M^K_\BC \ar[d]_{j_*}  &
\VMHS_F M^{K_1}_\BC \ar[d]^{(j_1)_*}  \\
\MHM_F M^K (\FS)_\BC \ar[r]^{\Ss} &
\MHM_F M^{K_1} (\Ses)_\BC \\}}
\]
which is compatible with that of \ref{2Cb}.
\end{Cor}

\begin{Proof}
By Theorem~\ref{2K}, we have
\[
(j_1)^* \Ss j_* \circ \mu_K = \mu_{K_1} \circ \Res^P_{P_1} \; .
\]
In order to see that the adjoint morphism
\[
\Ss j_* \circ \mu_K \longto (j_1)_* \circ \mu_{K_1} \circ \Res^P_{P_1}
\]
is an isomorphism, one may apply the (faithful)
forgetful functor to perverse sheaves on $M^{K_1} (\Ses)_\BC$. There, the
claim follows from Proposition~\ref{2Cb}.
\end{Proof}

We are ready to prove our main result:

\begin{Thm}\label{2H}
  There is a commutative diagram
\[
\vcenter{\xymatrix@R-10pt{ 
D^b \left( \Rep_F P \right) \ar[r]^{\Res^P_{P_1}} \ar[d]_{\mu_K} &
D^b \left( \Rep_F P_1 \right) \ar[r]^{R (\;)^{\langle \sigma \rangle}} &
D^b \left( \Rep_F P_{1,[\sigma]} \right) 
                        \ar[d]^{\mu_{\pi_{[\sigma]} (K_1)}} \\
D^b \left( \MHM_F M^K_{\BC} \right) \ar[rr]^{i^{\ast} j_{\ast}[-c]} &&
D^b \left( \MHM_F \MpC \right) \\}}
\]
which is compatible with that of \ref{2F}.
\end{Thm}

\begin{Proof}
According to \cite{Sa}~2.30, we have the equality
\[
i^* = i_1^* \circ \Ss \; .
\]
Together with Corollary~\ref{2Ka}, this reduces us to the case $P=P_1$.
Now observe that
$(i_1)_*$ and $(j_1)_*$ are exact functors on the level of abelian categories
$\MHM_F$ 
(\cite{BBD}~Cor.~4.1.3). $(i_1)^*$ is the left adjoint of $(i_1)_*$ 
on the level of
$D^b ( \MHM_F )$. It follows formally that the zeroeth cohomology functor
\[
\CH^0 (i_1)^*: \MHM_F M^{K_1} (\Ses)_\BC \longto \MHM_F \MpC 
\]
is right exact, and that $(\CH^0 (i_1)^*, (i_1)_*)$ constitutes an adjoint pair of
functors on the level of $\MHM_F$. In particular, there is an adjunction
morphism
\[
\id \longto (i_1)_* \CH^0 (i_1)^*
\]
of functors on $\MHM_F M^{K_1} (\Ses)_\BC$, 
which induces a morphism of functors on
\[
K^b \left( \MHM_F M^{K_1} (\Ses)_\BC \right) \; ,
\]
the homotopy category of complexes in $\MHM_F M^{K_1} (\Ses)_\BC$. 
Denote by $q$ the localization
functor from the homotopy to the derived category. We get a morphism in
\[
\begin{array}{cccc}
\Hom \left( q, q \circ (i_1)_* \CH^0 (i_1)^* \right) & = 
   & \Hom \left( q, (i_1)_* \circ q \circ \CH^0 (i_1)^* \right) & \\
& = & \Hom \left( (i_1)^* \circ q, q \circ \CH^0 (i_1)^* \right) & ,
\end{array}
\] 
where $\Hom$ refers to morphisms of exact functors. Composition with the
exact functor $(j_1)_* \circ \mu_{K_1}$ gives a morphism
\[
\eta' \in \Hom \left( (i_1)^* (j_1)_* \circ \mu_{K_1} \circ q, 
                      q \circ \CH^0 (i_1)^* (j_1)_* \circ \mu_{K_1} \right) \; .
\]
Assuming the existence of the {\it total left derived functor} 
\[
L \left( \CH^0 (i_1)^* (j_1)_* \circ \mu_{K_1} \right): 
       D^b \left( \Rep_F P_1 \right) \longto 
       D^b \left( \MHM_F \MpC \right)
\]
for a moment (see (a) below), its universal property (\cite{V}~II.2.1.2) says
that the above $\Hom$ equals
\[
\Hom \left( (i_1)^* (j_1)_* \circ \mu_{K_1}, 
            L \left( \CH^0 (i_1)^* (j_1)_* \circ \mu_{K_1} \right) \right) \; .
\]
Denote by
\[
\eta: (i_1)^* (j_1)_* \circ \mu_{K_1} \longto 
L \left( \CH^0 (i_1)^* (j_1)_* \circ \mu_{K_1} \right) 
\]
the morphism corresponding to $\eta'$.
It remains to establish the following claims:
\begin{itemize}
\item[(a)] The functor
\[
\CH^0 (i_1)^* (j_1)_* \circ \mu_{K_1}: \Rep_F P_1 \longto \MHM_F \MpC
\]
is left derivable.\\
\item[(b)] There is a canonical isomorphism between the total 
left derived functor
\[
L \left( \CH^0 (i_1)^* (j_1)_* \circ \mu_{K_1} \right)
\]
and
\[
\mu_{\pi_{[\sigma]} (K_1)} \circ R (\;)^{\langle \sigma \rangle} [c] \; .
\]
\item[(c)] $\eta$ is an isomorphism.
\end{itemize}

For (a) and (b), observe that up to a twist by $c$, 
the variation 
\[
\CH^0 (i_1)^* (j_1)_* \circ \mu_K (\BV)
\] 
on $\MpC$ is given by the co-invariants of $\BV$ under the local monodromy.
This is a general fact about the degeneration of variations along a divisor
with normal crossings; see e.g.\ the discussion preceding \cite{HZ}~(4.4.8).
By \cite{K}~Thm.~6.10, up to a twist
by $c$ (corresponding to the highest exterior power of $\Lie \langle \sigma \rangle$), the co-invariants are identical to $H^c (\langle \sigma \rangle, \;\;)$.

We are thus reduced to showing that 
the functor $H^c (\langle \sigma \rangle, \;\;)$ is left derivable, with
total left derived functor
$R (\;)^{\langle \sigma \rangle} [c]$. But this follows from standard facts
about Lie algebra homology and cohomology (see e.g.\ \cite{K}~Thm.~6.10 and
its proof).

(c) can be shown after applying the 
forgetful functor to perverse sheaves. There, the
claim follows from \ref{2F}.
\end{Proof}

\begin{Rem}\label{2I}
If $(P, \FX)$ is pure, and $c = \dim \langle \sigma \rangle$ is maximal, 
i.e., equal to $\dim U_1$, then Theorem~\ref{2H} is
equivalent to \cite{HZ}~Thm.~(4.4.18). In fact, 
by \ref{2H}, the recipe to compute $H^q i^{\ast} j_{\ast} 
\circ \mu_K (\BV)$ given on pp.~286/287 of \cite{HZ} generalizes as follows: The complex
\[
C^{\bullet} = \Lambda^{\bullet} (\Lie \langle \sigma \rangle )^{\ast} 
\otimes_F \BV
\]
carries the diagonal action of $P_1$ 
(where the action on $\Lie \langle \sigma \rangle$ is via conjugation).
The induced action on the cohomology objects $H^q C^{\bullet}$ factors through $P_{1,[\sigma]}$ and gives the right Hodge structures via $\mu_{\pi_{[\sigma]} (K_1)}$. In \cite{HZ}, the Hodge and weight filtrations on $C^{\bullet}$ corresponding to the action of $P_1$ are made explicit.
\end{Rem}

\begin{Rem}\label{2L}
Because of \ref{1O}~(b), the isomorphism of Theorem \ref{2H} does not depend on the cone decomposition $\FS$, which contains $\sigma \times \{ p \}$. We leave it to the reader to formulate and prove
results like \cite{P2}~(4.8.5) on the behaviour of the isomorphism of
\ref{2H} under change of the group $K$, and of the element $p$.
\end{Rem}

Let us conclude the section with a statement on transitivity of degeneration.
In addition to the data used so far, fix a face $\tau$ of $\sigma$. Write
\[
i_{\tau}: \Mpt \into M^K(\FS) \; .
\]
$\Mps$ lies in the closure of $\Mpt$ inside $M^K(\FS)$. Adjunction gives
a morphism
\[
i^* j_* \circ \mu_K \longto i^* (i_{\tau})_* (i_{\tau})^* j_* \circ \mu_K
\]
of exact functors from $D^b (\Rep_F P )$ to 
$D^b \left( \MHM_F \MpC \right)$.

\begin{Prop}\label{2M}
This morphism is an isomorphism. 
\end{Prop}

\begin{Proof}
This can be checked on the level of local systems. There, it follows from 
Theorem~\ref{1K}~(i), and standard facts about degenerations 
along strata in torus embeddings.
\end{Proof}

%
%

\section{Higher direct images for $\ell$-adic sheaves}
\label{3}


The main result of this section (Theorem~\ref{3J}) provides an $\ell$-adic analogue of the formula of \ref{2H}. The main ingredients of the proof are the machinery developed in \cite{P2}, and our knowledge of the local situation (\ref{1M}). \ref{3E}--\ref{3Fa} are concerned with the problem of extending certain infinite families of \'etale sheaves to ``good'' models of a Shimura variety. We conclude by discussing mixedness of the $\ell$-adic sheaves obtained via the canonical construction.\\

With the exception of condition $(\cong)$, which will not be needed, we fix the same set of geometric data as in the beginning of Section \ref{2}.
In particular, the cone $\sigma$ is assumed smooth, the group $K$ is neat, and $(P, \FX)$ satisfies condition $(+)$.\\

Define $\tilde{M} (\FS)$ as the inverse limit of all
\[
M^{K'} (\FS) = M^{K'} (P , \FX , \FS)
\]
for open compact $K' \le K$. The group $K$ acts on $\tilde{M} (\FS)$, and
\[
M^K (\FS) = \tilde{M} (\FS) / K \; .
\]

Inside $\tilde{M} (\FS)$ we have the inverse limit $\tilde{M}$ of
\[
M^{K'} = M^{K'} (P , \FX) \; , \quad K' \le K \; ,
\]
and the inverse limit $\tilde{M}_{[\sigma]}$ of all
\[
M^{K'_{1,[\sigma]}} = M^{K'_{1,[\sigma]}} (\Pes , \Xes)
\]
for open compact $K'_{1,[\sigma]} \le K_{1,[\sigma]} := \pi_{[\sigma]} (K_1)$.
We get a commutative diagram
\begin{myequation}\label{3eq}
\vcenter{\xymatrix@R-10pt{ 
\tilde{M} \ar@{^{ (}->}[r]^{\tilde{j}} \ar[d] &
\tilde{M} (\FS) \ar@{<-^{ )}}[r]^{\tilde{i}} \ar[d] &
\tilde{M}_{[\sigma]} \ar[d] \\
\tilde{M} / K_1  \ar@{^{ (}->}[r]^{j'} \ar[d]_{\varphi} &
\tilde{M} (\FS) / K_1 \ar@{<-^{ )}}[r]^{i'} \ar[d]_{\tilde{\varphi}} &
\tilde{M}_{[\sigma]} / K_{1 , [\sigma]} \ar@{=}[d] \\
M^K = \tilde{M} / K  \ar@{^{ (}->}[r]^{\!\!\!\!\!\!\!\! j} &
M^K (\FS) = \tilde{M} (\FS) / K \ar@{<-^{ )}}[r]^{\quad\quad i} &
\Mps \\}}
\end{myequation}

\begin{Prop}\label{3B}
  The morphism
\[
\tilde{\varphi} : \tilde{M} (\FS) / K_1 \longrightarrow M^K (\FS) = \tilde{M} (\FS) / K
\]
is \'etale near the stratum
\[
\Mps = \tilde{M}_{[\sigma]} / K_{1,[\sigma]} \; .
\]
\end{Prop}

\begin{Proof}
  By Theorem~\ref{1M}, the map $\tilde{\varphi}$ induces an isomorphism of the respective formal completions along our stratum. The claim thus follows from \cite{EGAIV4}~Prop.~(17.6.3).
\end{Proof}

Let $\Tor \Mod_K$ be the category of all continuous discrete torsion $K$-modules. The left vertical arrow of (\ref{3eq})
gives an evident functor
\[
\mu_K : \Tor \Mod_K \longrightarrow \Et M^K
\]
into the category of \'etale sheaves on $M^K$; since $K$ is neat, this functor is actually an exact tensor functor with values in the category of lisse sheaves. Similar remarks apply to $K_1$ or $\pi_{[\sigma]} (K_1)$ in place of $K$. We are interested in the behaviour of the functor
\[
i^{\ast} j_{\ast} : D^+ (\Et M^K ) \longrightarrow D^+ ( \Et \Mps )
\]
on the image of $\mu_K$. 
From \ref{3B}, we conclude:

\begin{Prop}\label{3C}
  (i) The two functors
\[
i^{\ast} j_{\ast} \quad , \quad (i')^{\ast} j'_{\ast} \circ \varphi^{\ast}:
D^+ (\Et M^K) \longrightarrow D^+ ( \Et \Mps )
\]
are canonically isomorphic.\\
(ii) The two functors
\[
i^{\ast} j_{\ast} \circ \mu_K \; , \; (i')^{\ast} j'_{\ast} \circ \mu_{K_1} \circ \Res^K_{K_1}: 
D^+ (\Tor \Mod_K) \longrightarrow D^+ ( \Et \Mps ) 
\]
are canonically isomorphic. Here, $\Res^K_{K_1}$ denotes the pullback via the monomorphism
\[
K_1 \longrightarrow K \; , \; k_1 \longmapsto p^{-1} \cdot k_1 \cdot p \; .
\]

\end{Prop}

\begin{Proof}
  (i) is smooth base change, and (ii) follows from (i). 
\end{Proof}

Write $K_\sigma$ for $\ker (\pi_{[\sigma]} \, |_{K_1}) = K_1 \cap \langle \sigma \rangle (\BA_f)$.

\begin{Thm}\label{3A}
  There is a commutative diagram 
\[
\vcenter{\xymatrix@R-10pt{ 
D^+ \left( \Tor \Mod_K \right) \ar[r]^{\Res^K_{K_1}} \ar[d]_{\mu_K} &
D^+ \left( \Tor \Mod_{K_1} \right) 
                      \ar[r]^{R (\;)^{\! K_\sigma}} &
D^+ ( \Tor \Mod_{\pi_{[\sigma]} (K_1)} )
                      \ar[d]^{\mu_{\pi_{[\sigma]} (K_1)}} \\
D^+ \left( \Et M^K \right) \ar[rr]^{i^{\ast} j_{\ast}} &&
D^+ \left( \Et \Mps \right) \\}}
\]
Here, $R (\;)^{K_\sigma}$ refers to continuous group cohomology of $K_\sigma$.
\end{Thm}

\begin{Proof}
We need to show that the diagram
\[
\vcenter{\xymatrix@R-10pt{ 
D^+ \left( \Tor \Mod_{K_1} \right) 
           \ar[r]^{R (\;)^{K_\sigma}} 
           \ar[d]_{\mu_{K_1}} &
D^+ ( \Tor \Mod_{\pi_{[\sigma]}(K_1)} ) 
                          \ar[d]^{\mu_{\pi_{[\sigma]} (K_1)}} \\
D^+ ( \Et \tilde{M} / K_1 )
           \ar[r]^{(i')^{\ast} j'_{\ast}} &
D^+ \left( \Et \Mps \right) \\}}
\]
commutes. The proof of this statement makes use of the full machinery developed in the first two section of \cite{P2}.

In fact, \cite{P2}~Prop.~(4.4.3) contains the analogous statement for the (coarser) stratification of $M^K (\FS)$ induced from the canonical stratification of the {\it Baily--Borel compactification} of $M^K$. One faithfully imitates the proof, observing that \cite{P2}~(1.9.1) can be applied because the upper half of (\ref{3eq})
is cartesian up to nilpotent elements. The statement on ramification along a stratum in \cite{P2}~(3.11) holds for arbitrary, not just pure Shimura data.
\end{Proof}

\begin{Rem}\label{3D}
 Because of Remark~\ref{1O}~(b), the isomorphism of \ref{3A} does not depend on the cone decomposition $\FS$ containing $\sigma \times \{ p \}$.
\end{Rem}

Fix a set $\CT \subset \Tor \Mod_K$, let $E = E (P, \FX)$ be the field of definition of our varieties, and write $O_E$ for its ring of integers. Consider a {\it model}
\[
\CM^K \stackrel{j}{\hookrightarrow} \CM^K (\FS) \stackrel{i}{\hookleftarrow} \CM^{\pi_{[\sigma]} (K_1)}
\]
of
\[
M^K \stackrel{j}{\hookrightarrow} M^K (\FS) \stackrel{i}{\hookleftarrow} \Mps
\]
over $O_E$, i.e., normal schemes of finite type over $O_E$, an open immersion $j$ and an immersion $i$ whose generic fibres give the old situation over $E$; we require also that the generic fibres lie dense in their models.
(Finitely many special fibres of our models might be empty.)

Assume
\begin{enumerate}
\item [(1)] All sheaves in $\mu_K (\CT)$ extend to lisse sheaves on $\CM^K$.
\item [(2)] For any $S \in \mu_K (\CT)$ and any $q \ge 0$, the extended sheaf $\CS$ on $\CM^K$ satisfies the following:
\[
i^{\ast} R^q j_{\ast} \CS \in \Et \CM^{\pi_{[\sigma]} (K_1)} \; \mbox{is lisse.}
\]
\end{enumerate}
Then the generic fibre of $i^{\ast} R^q j_{\ast} \CS$ is necessarily equal to $i^{\ast} R^q j_{\ast} S$, i.e., it is given by the formula of \ref{3A}. So $i^{\ast} R^q j_{\ast} \CS$ is the unique lisse extension of $i^{\ast} R^q j_{\ast} S$ to $\CM^{\pi_{[\sigma]} (K_1)}$. Observe that if $\CT$ is finite, then conditions (1) and (2) hold after passing to an open sub-model of any given model.\\
If $\CT$ is an abelian subcategory of $\Tor \Mod_K$ and (1) holds, then (2) needs to be checked only for the simple noetherian objects in $\CT$.\\

Let us show how to obtain a model as above for a \emph{particular} choice of $\CT$:\\
Fix a prime $\ell$, write
\[
\pr_\ell : P (\BA_f) \longrightarrow P (\BQ_\ell)
\]
and $K_\ell := \pr_\ell (K)$. Denote by $\CT_\ell \subset \Tor \Mod_{K_\ell} \subset \Tor \Mod_K$ the abelian subcategory of $\BZ_\ell$-torsion $K_\ell$-modules. The quotient $K_\ell$ of $K$ corresponds to a certain part of the ``Shimura tower''
\[
( M^{K'} )_{K'} \; ,
\]
namely the one indexed by the open compact $K' \le K$ containing the kernel of $\pr_\ell \, |_K$. According to \cite{P2}~(4.9.1), the following is known:

\begin{Prop}\label{3E}
  There exists a model $\CM^K$ such that all the sheaves in
\[
\mu_K (\CT_\ell)
\]
extend to lisse sheaves on $\CM^K$. Equivalently, the whole \'etale $K_\ell$-covering of $M^K$ considered above extends to an \'etale $K_\ell$-covering of $\CM^K$.
\end{Prop}

\begin{Proof}
Write $L$ for the product of $\ell$ and the primes dividing the order of the
torsion elements in $K_\ell$; thus $K_\ell$ is a pro-$L$-group.
Let $S$ be a finite set in $\Spec O_E$ containing the prime factors of $L$, and $\CM^K$ a model of $M^K$ over $O_S$ which is the complement of an $NC$-divisor relative to $O_S$ in a smooth, proper $O_S$-scheme.

We give a construction of a suitable enlargement $S'$ of $S$ such that the claim holds for the restriction of $\CM^K$ to $O_{S'}$.\\

First, assume that $P$ is a torus. Recall (\cite{P1}~2.6) that the Shimura varieties associated to tori are
finite over their reflex field. 
Since Shimura varieties are normal, each $M^K$ is the spectrum of a 
finite product $E_K$ of number fields. But then the $K_\ell$-covering corresponds to an {\it abelian} $K_\ell$-extension
\[
\tilde{E} / E_K \; .
\]
By looking at the kernel of the reduction map to $\GL_N (\BZ / \ell^f \BZ)$, $\ell^f \ge 3$, one sees that there is an intermediate extension
\[
\tilde{E} / F / E_K
\]
finite over $E_K$, such that $\tilde{E} / F$ is a $\BZ^r_\ell$-extension. Hence the only primes that ramify in $\tilde{E} / F$ are those over $\ell$, and one adds to $S$ the finitely many primes which ramify in $F / E_K$.\\

In the general case, choose an embedding
\[
e : (T , \CY) \longrightarrow (P , \FX)
\]
of Shimura data, with a torus $T$ 
such that $E = E(P,\FX)$ is contained in $E(T, \CY)$
(\cite{P1}~Lemma~11.6), and finitely many $K^T_m \le T (\BA_f)$ and $p_m \in P (\BA_f)$ such that the maps
\[
[\cdot p_m] \circ [e] : M^{K^T_m} (T , \CY) \longrightarrow M^K (P , \FX) 
\]
are defined and meet all components of $M^K$ (\cite{P1}~Lemma~11.7).
Each $M^{K^T_m}$ equals the spectrum of a product $F_m$ of number fields.

Define $x_m \in M^K (F_m)$ as the image of $[\cdot p_m] \circ [e]$. 
Let $S_m \subset \Spec O_{F_m}$ denote the set of bad primes for $M^{K^T_m}$ and $(K^T_m)_\ell$, plus a suitable finite set such that $x_m$ extends to a section of $\CM^K$ over $O_{S_m}$. 

Enlarge $S = S((T, \CY),e,p_m)$ so as to contain all primes which ramify in some $F_m$, and those below a prime in some $S_m$. We continue to write $S$ for the
enlargement, and $\CM^K$ and $x_m$ for the objects obtained via restriction to $O_S$.\\

We claim that with these choices, the whole \'etale $K_\ell$-covering of 
$M^K$ extends to an \'etale $K_\ell$-covering of $\CM^K$.

Let $M^0$ and $\CM^0$ be connected components of $M^K$ and $\CM^K$. We have to show that the map
\[
s : \pi_1 (M^0) \longrightarrow K_\ell
\]
given by the $K_\ell$-covering factors through the epimorphism
\[
\beta : \pi_1 (M^0) \longonto \pi_1 (\CM^0) \; .
\]
There is an $m$ and intermediate field extensions
\[
F_m / F' / F / E
\]
such that $M^0$ is a scheme over $F$ with geometrically connected fibres,
and such that $x_m$ induces an $F'$-valued point of $M^0$.
Since $\CM^0$ is normal, $\CM^0$ is a scheme over the integral closure $O_{S_F}$ of $O_S$. By \cite{SGA1XIII}~4.2--4.4, there is a commutative diagram of exact sequences
\[
\vcenter{\xymatrix@R-10pt{ 
1 \ar[r] & 
\pi_1 (\overline{M^0}) \ar[r] \ar[d]_{\alpha} &
\pi_1 (M^0) \ar[r] \ar[d]_{\beta'} &
\Gal (\overline{F} / F) \ar[r] \ar[d]^{\gamma} &
1 \\
1 \ar[r] & 
\pi^L_1 (\overline{M^0}) \ar[r] &
\pi'_1 (\CM^0) \ar[r] &
\pi_1 (\Spec O_{S_F}) \ar[r] &
1 \\}}
\]
Here, $\pi^L_1 (\overline{M^0})$ is the largest pro-$L$-quotient of $\pi_1 (\overline{M^0})$, the fundamental group of $\overline{M^0} := M^0 \otimes_F \overline{F}$, and $\pi'_1 (\CM^0)$ is a suitable quotient of $\pi_1 (\CM^0)$.
Hence all vertical arrows are surjections.

Clearly $\ker \alpha$ is contained in $\ker s$; we thus get a map
\[
s' : \pi_1 (M^0) / \ker \alpha \longrightarrow K_\ell \; .
\]
We have to check that
\[
\ker \gamma = \ker \beta' / \ker \alpha \subset \pi_1 (M^0) / \ker \alpha
\]
is contained in $\ker s'$. But $\ker \gamma$ remains unchanged under passing to the extension $F' / F$, which is unramified outside $S_F$. There, the corresponding exact sequence splits thanks to the existence of $x_m$.\\
The map
\[
\ker \gamma \longrightarrow \pi_1 (M^0)
\]
is induced by pullback via $[\cdot p_m] \circ [e]$, and by construction its image is contained in $\ker s$.
\end{Proof}

This takes care of condition (1).

\begin{Lem}\label{3F}
Up to isomorphism, there are only finitely many simple objects in $\CT_\ell$.
\end{Lem}

\begin{Proof}
There is a normal subgroup $K'_\ell \le K_\ell$ of finite index which is a projective limit of $\ell$-groups. Write
$\CT'_\ell$ for the subcategory of $\Tor \Mod_{K'_\ell}$ of $\BZ_\ell$-torsion modules. Since any element of order $\ell^n$ in $\GL_r (\BF_\ell)$ is unipotent,
any simple non-trivial object in $\CT'_\ell$ is isomorphic to the trivial representation $\BZ / \ell \BZ$ of $K'_\ell$.

Therefore, the simple objects in $\CT_\ell$ all occur in the Jordan--H\"older decomposition of
\[
\Ind^{K_\ell}_{K'_\ell} \Res^{K_\ell}_{K'_\ell} (\BZ / \ell\BZ) \; .
\]
\end{Proof}

\begin{Prop}\label{3Fa}
Conditions (1) and (2) hold for 
a suitable open sub-model of
any model as in \ref{3E}.
\end{Prop}

\begin{Proof}
By generic base change (\cite{SGAhalbTh}~Thm.~1.9), condition (2) can be
achieved for any single constructible sheaf $\CS$ on $\CM^K$,
which is lisse on the generic fibre. 
The claim
follows from \ref{3F} by applying the long exact sequences associated
to $i^\ast R j_\ast$.
\end{Proof}

Fix a finite extension $F = F_{\lambda}$ of $\BQ_\ell$. By passing to projective limits, we get an exact tensor functor
\[
\mu_K : \Rep_F P \longrightarrow \Et^l_F M^K
\]
into the category of lisse $\lambda$-adic sheaves on $M^K$ (\cite{P2}~(5.1)). 
We refer to $\mu_K$ as the {\it canonical construction} of $\lambda$-adic
sheaves from representations of $P$.
Denote by $D^b_\con(?,F)$ Ekedahl's bounded ``derived'' category of constructible
$F$-sheaves (\cite{E}~Thm.~6.3). Consider the functor
\[
i_{\ast} j^{\ast} : D^b_\con (M^K,F) \longrightarrow D^b_\con ( \Mps ,F) \; .
\]
From Theorem~\ref{3A}, we obtain the main result of this section:

\begin{Thm}\label{3J}
  There is a commutative diagram
\[
\vcenter{\xymatrix@R-10pt{ 
D^b \left( \Rep_F P \right) \ar[r]^{\Res^P_{P_1}} \ar[d]_{\mu_K} &
D^b \left( \Rep_F P_1 \right) \ar[r]^{R (\;)^{\langle \sigma \rangle}} &
D^b \left( \Rep_F P_{1,[\sigma]} \right) 
                        \ar[d]^{\mu_{\pi_{[\sigma]} (K_1)}} \\
D^b_\con \left( M^K, F \right) \ar[rr]^{i^{\ast} j_{\ast}} &&
D^b_\con \left( \Mps, F \right) \\}}
\]
Here, $\Res^P_{P_1}$ denotes the pullback via the monomorphism
\[
P_{1,F} \longrightarrow P_F \; , \; p_1 \longmapsto \pi_\ell (p)^{-1} \cdot p_1 \cdot \pi_\ell (p) \; ,
\]
and $R (\;)^{\langle \sigma \rangle}$ is Hochschild cohomology of the unipotent group $\langle \sigma \rangle$.
\end{Thm}

\begin{Proof}
  Since $\langle \sigma \rangle$ is unipotent, $R (\;)^{K_\sigma}$ and $R (\;)^{\langle \sigma \rangle}$ agree. 
\end{Proof}

Let us note a refinement of the above. Consider smooth models
\[
\CM^K \stackrel{j}{\hookrightarrow} \CM^K (\FS) \stackrel{i}{\hookleftarrow} \CM^{\pi_{[\sigma]} (K_1)}
\]
satisfying conditions (1), (2) for $\CT_\ell$. Thus 
all the sheaves in the image of $\mu_K$ extend to $\CM^K$; in particular they can be considered as (locally constant)
{\it perverse $F$-sheaves} in the sense of \cite{Hu2}:
\[
\mu_K : \Rep_F P \longrightarrow \Perv_F \CM^K \subset D^b_\con (\FU \CM^K , F) 
\]
(notation as in \cite{Hu2}).
Consider the functor
\[
i_{\ast} j^{\ast} : D^b_\con (\FU \CM^K , F) \longrightarrow 
D^b_\con ( \FU \CM^{\pi_{[\sigma]} (K_1)} , F ) \; .
\]

\begin{Var}\label{3K}
  There is a commutative diagram
\[
\vcenter{\xymatrix@R-10pt{ 
D^b \left( \Rep_F P \right) \ar[r]^{\Res^P_{P_1}} \ar[d]_{\mu_K} &
D^b \left( \Rep_F P_1 \right) \ar[r]^{R (\;)^{\langle \sigma \rangle}} &
D^b \left( \Rep_F P_{1,[\sigma]} \right) 
                        \ar[d]^{\mu_{\pi_{[\sigma]} (K_1)}} \\
D^b_\con \left(\FU \CM^K , F \right) \ar[rr]^{i^{\ast} j_{\ast}[-c]} &&
D^b_\con \left(\FU \CM^{\pi_{[\sigma]} (K_1)} , F \right) \\}}
\]
\end{Var}

\begin{Rem}\label{3L}
  As in \ref{3A}, the isomorphism 
\[
\mu_{\pi_{[\sigma]} (K_1)} \circ 
R (\;)^{\langle \sigma \rangle} \circ \Res^P_{P_1}   \isoto
i^{\ast} j_{\ast} \circ \mu_K [-c]
\]
does not depend on the cone decomposition $\FS$ containing $\sigma \times \{ p \}$. It is possible, as in \cite{P2}~(4.8.5), to identify the effect on
the isomorphism of change of the group $K$ and of the element $p$. Similarly,
one has an $\ell$-adic analogue of Proposition~\ref{2M}.
\end{Rem}

In the above situation, consider the {\it horizontal stratifications} (\cite{Hu2}~page~110) ${\bf S} = \{ \CM^K \}$ of $\CM^K$ and ${\bf T} = \{ \CM^{\pi_{[\sigma]} (K_1)} \}$ of $\CM^{\pi_{[\sigma]} (K_1)}$. Write $L_{\bf S}$ and $L_{\bf T}$ for the sets of extensions to the models of irreducible objects of $\mu_K (\Rep_F P)$ and $\mu_{\pi_{[\sigma]} (K_1)} (\Rep_F \Pes)$ respectively. In the terminology of \cite{Hu2}~Def.~2.8, we have the following:

\begin{Prop}\label{3M}
  $i_{\ast} j^{\ast}$ is $({\bf S} , L_{\bf S})$-to-$({\bf T} , L_{\bf T})$-admissible.
\end{Prop}

\begin{Proof}
  This is \cite{Hu2}~Lemma~2.9, together with Theorem \ref{3J}.
\end{Proof}

It is conjectured (\cite{LR}~\S\,6; \cite{P2}~(5.4.1); \cite{W}~4.2) that the image of $\mu_K$ consists of {\it mixed sheaves with a weight filtration}; furthermore, the filtration should be the one induced from the weight filtration of representations of $P$. Let us refer to this as the {\it mixedness conjecture} for $(P,\FX)$; cmp.\ \cite{P2}~(5.5)--(5.6) and \cite{W}~pp~112--116
for a discussion. The conjecture is known if every $\BQ$-simple factor of $G^{\ad}$ is {\it of abelian type} (\cite{P2}~Prop.~(5.6.2), \cite{W}~Thm.~4.6~(a)).

\begin{Prop}\label{3N}
  If the mixedness conjecture holds for $(P,\FX)$, 
then it holds for any rational boundary component $(P_1, \FX_1)$.
\end{Prop}

\begin{Proof}
  By \cite{W}~Thm.~4.6, it suffices to check that $\mu_{\pi_{[\sigma]} (K_1)} (\BW)$ is mixed for some faithful representation $\BW$ of $\Pes$. By \cite{Hm}~Thm.~11.2, there is a representation $\BV$ of $P$ and a one-dimensional subspace $\BV' \subset \BV$ such that
\[
\langle \sigma \rangle = \Stab_P (\BV') \; .
\]
Since $\langle \sigma \rangle$ is unipotent, we have
\[
\BV' \subset \BW := H^0 (\langle \sigma \rangle , \BV) \; .
\]
$\BW$ is a faithful representation of $\Pes$, and by Theorem~\ref{3J}, $\mu_{\pi_{[\sigma]} (K_1)} (\BW)$ is a cohomology object of the complex
\[
i^{\ast} j_{\ast} \circ \mu_K (\BV) \; .
\]
It is therefore mixed (\cite{D}~Cor.~6.1.11).
\end{Proof}

%
%

\end{document}